\documentclass{amsart}
\usepackage{amssymb,amsmath,amsthm,cite,verbatim,IEEEtrantools,url, macros}
\title{Large Subsets of Local Fields Not Containing Configurations}
\author{Robert Fraser}
\address{Department of Mathematics \\
University of British Columbia \\ Canada}
\email{rgf@math.ubc.ca}
\begin{document}
\begin{abstract}
For certain families of functions $\{f_{\ell}\}$ mapping $K^{nv_{\ell}} \to K^m$, where $K$ is a complete, nonarchimedean local field, we find a set $E$ of large Hausdorff dimension with the property that $f_{\ell}(x_1, \ldots, x_{v_{\ell}})$ is nonzero for any distinct points $x_1, \ldots, x_{v_{\ell}} \in E$. In particular, this result can be applied to show that the ring of integers of any local field contains a subset of Hausdorff dimension $1$ not containing any nondegenerate 3-term arithmetic progressions.
\end{abstract}
\maketitle
\section{Introduction}
We are interested in questions of the following form: given a family of functions $\{f_{\ell}\} : R^{nv_{\ell}} \to K^m$, where $K$ is a complete, nonarchimedean local field with ring of integers $R$, how large can the Hausdorff dimension of the set $E \subset R^n$ be if there do not exist any $v_{\ell}$-tuples $(x_1, \ldots, x_{v_{\ell}})$ with $x_1, \ldots, x_{v_{\ell}}$ distinct elements of $E$ such that $f_{\ell}(x_1, \ldots, x_{v_{\ell}}) = 0$ for any $\ell$? Versions of this question have been considered for the real numbers by Keleti \cite{K98, K08}, Maga \cite{M10}, M\'ath\'e \cite{M17}, and Pramanik and the author \cite{FP16}. The objective of this paper is to establish local field versions of the results \cite{M17} and \cite{FP16}, which will entail results similar to those occurring in \cite{K98} and \cite{K08} as a special case.

In order to draw an analogy between the theorems of this paper and their Euclidean counterparts, we quote the Euclidean versions verbatim, up to notational differences. We begin by quoting a theorem of M\'ath\'e \cite{M17}:

\begin{quote}
\begin{mythm*}[M\'ath\'e]
Let $n \geq 1$. Let $L$ be a countable set. For each $\ell \in L$, let $v_{\ell}$ be a positive integer, and let $P_{\ell} : \mathbb{R}^{n v_{\ell}} \to \mathbb{R}$ be a (non identically zero) polynomial in $n v_{\ell}$ variables with rational coefficients. Assume that $d$ is the maximum degree of the polynomials $P_{\ell} (\ell \in L)$. Then there exists a compact set $E \subset \mathbb{R}^n$ of Hausdorff dimension $n/d$ such that for every $j \in J$, $E$ does not contain $v_{\ell}$ distinct points $x_1, \ldots, x_{v_{\ell}}$ satisfying $P_j(x_1, \ldots, x_{v_{\ell}}) = 0$. 
\end{mythm*}
\end{quote}

We now state the analogous result for local fields. The \emph{height} appearing in the theorem statement is defined in Section 2; in the case of $R = \mathbb{Z}_p$, elements of finite height are ordinary integers.
\begin{mythm}\label{mathe}
Let $\{f_{\ell}\} : R^{nv_{\ell}} \to R$ be a countable family of nonzero polynomials of degree at most $d$ whose coefficients are elements of $R$ of finite height. Suppose further that either $R$ has characteristic $0$ or $d < \text{char }R$. Then there exists a set $E \subset R$ of Hausdorff dimension $\frac{n}{d}$ and Minkowski dimension $n$ such that, for all ${\ell}$, the set $E$ does not contain $v_{\ell}$ distinct points $x_1, \ldots, x_{v_{\ell}}$ such that $f_{\ell}(x_1, \ldots, x_{v_{\ell}}) = 0$. 
\end{mythm}
This theorem can be applied to the function $x_1 - 2x_2 + x_3$ to show the existence of a subset of Hausdorff dimension $1$ of any local field that does not contain any $3$-term arithmetic progressions with distinct elements. This contrasts with a recent result of Ellenberg and Gijswijt \cite{EG17}, which states that for finite dimensional vector spaces over (say) $\mathbb{Z}/3 \mathbb{Z}$ of dimension $N$, there is a constant $C < 3$ such that any set of at least $C^N$ elements contains a $3$-term arithmetic progression. In fact, there is a subset $E$ of any complete, nonarchimedean local field $K$ of Hausdorff dimension $1$ with no $v$-tuple of distinct points satisfying any linear equation of the form $a_1x_1 + \cdots + a_v x_v = 0$ where $a_1, \ldots, a_v \subset R$ have finite height. This implies in particular that if the characteristic of the field $K$ is equal to zero, then there are no solutions to any equations of this form with integer coefficients; that is, $E$ is linearly independent over $\mathbb{Q}$.

The assumption on the degree $d$ in the case where $R$ is a ring of finite characteristic is important because the argument requires some higher-order partial derivative of each polynomial to be nonvanishing. This may not happen in finite characteristic for polynomials of degree greater than or equal to the characteristic.

Next, we quote two theorems verbatim from \cite{FP16}. In these theorems, the number of variables $v$ is fixed, but the functions need not be polynomials.
\begin{quote}
\begin{mythm*}[Fraser-Pramanik]
For any $\eta > 0$ and integer $v \geq 3$, let $f_{\ell}: \mathbb{R}^v \to \mathbb{R}$ be a countable family of functions in $v$ variables with the following properties: 
\begin{enumerate}[(a)]
\item There exists $r_{\ell} < \infty$ such that $f_{\ell} \in C^{r_{\ell}}([0, \eta]^v)$,  
\vskip0.1in 
\item For each $\ell$, some partial derivative of $f_{\ell}$ of order $r_{\ell} \geq 1$ does not vanish at any point of $[0, \eta]^v$. 
\end{enumerate} 
\vskip0.1in 
\noindent Then there exists a set $E \subseteq [0, \eta]$ of Hausdorff dimension at least $\frac{1}{v-1}$ and Minkowski dimension 1 such that $f_{\ell}(x_1, \ldots, x_v)$ is not equal to zero for any $v$-tuple of distinct points $x_1, \ldots, x_v \in E$ and any function $f_{\ell}$.
\end{mythm*}
\end{quote}
This theorem has a multidimensional analogue, but an additional assumption is needed on the derivative. This multidimensional theorem is quoted below.
\begin{quote}
\begin{mythm*}[Fraser-Pramanik]
Fix $\eta > 0$ and positive integers $m,n,v$ such that $v \geq 3$, and $m \leq n(v-1)$. Let $f_{\ell} : (\mathbb{R}^{n})^v \to \mathbb{R}^m$ be a countable family of $C^{2}$ functions with the following property: the derivative $Df_{\ell}(x_1, \ldots, x_v)$ has full rank at every point $(x_1, \cdots, x_v)$ in the zero set of $f_{\ell}$ such that $x_r \ne x_s$ for all $r \ne s$. 
\vskip0.1in
\noindent Then there exists a set $E \subseteq [0, \eta]^{n}$ of Hausdorff dimension at least $\frac{m}{v-1}$ and Minkowski dimension $n$ such that $f_{{\ell}}(x_1, \ldots, x_v)$ is not equal to zero for any $v$-tuple of distinct points $x_1, \ldots, x_v \in E^n$ and any function $f_{{\ell}}$.
\end{mythm*}
\end{quote}
We will establish the following local-field analogues of these two theorems:
\begin{mythm}\label{non-simul-thm}
For any ball $B \subset R$ and integer $v \geq 3$, let $f_{{\ell}}: R \to K$ be a countable family of functions in $v$ variables with the following properties: 
\begin{enumerate}[(a)]
\item There exists $r_{\ell} < \infty$ such that $f_{\ell}$ is $r_{\ell}$ times strictly differentiable.
\vskip0.1in 
\item For each ${\ell}$, some partial derivative of $f_{\ell}$ of order $r_{{\ell}} \geq 1$ does not vanish at any point of $B^v$. 
\end{enumerate} 
\vskip0.1in 
\noindent Then there exists a set $E \subseteq B$ of Hausdorff dimension at least $\frac{1}{v-1}$ and Minkowski dimension 1 such that $f_{{\ell}}(x_1, \ldots, x_v)$ is not equal to zero for any $v$-tuple of distinct points $x_1, \ldots, x_v \in E$ and any function $f_{{\ell}}$.
\end{mythm}
\begin{mythm}\label{non-simul-thm vector-valued}
Fix a ball $B \subset R$ and positive integers $m,n,v$ such that $v \geq 3$, and $m \leq n(v-1)$. Let $f_{{\ell}} : (R^{n})^v \to K^m$ be a countable family of twice strictly differentiable functions with the following property: the first derivative $Df_{\ell}$ has full rank on the portion of the zero set of $f_{\ell}$ contained in $B^{nv}$ for every ${\ell}$. 
\vskip0.1in
\noindent Then there exists a set $E \subseteq B^{n}$ of Hausdorff dimension at least $\frac{m}{v-1}$ and Minkowski dimension $n$ such that $f_{{\ell}}(x_1, \ldots, x_v)$ is not equal to zero for any $v$-tuple of distinct points $x_1, \ldots, x_v \in E^n$ and any function $f_{{\ell}}$.
\end{mythm}
Finally, we quote one last theorem from \cite{FP16}. This theorem discusses uncountable families of configurations with the same linearization at each diagonal point.
\begin{quote}
\begin{mythm*}[Fraser-Pramanik]
Given any constant $C > 0$ and a vector $\alpha \in \mathbb R^v$ such that 
\[\alpha \cdot u \ne 0  \text{ for every } u \in \{ 0, 1\}^v \text{ with } u \ne 0,  u \neq (1, 1, \ldots, 1)\]
and such that
\[\sum_{j=1}^v \alpha_j = 0,\]
there exists a positive constant $c(\alpha)$ and a set $E = E(C, \alpha) \subseteq [0,1]$ of Hausdorff dimension $c(\alpha) > 0$ with the following property. 
\vskip0.1in
\noindent The set $E$ does not contain any nontrivial solution of the equation 
\[f(x_1, \cdots, x_v) = 0, \qquad (x_1, \cdots, x_v) \text{ not all identical}, \] 
for any $C^2$ function $f$ of the form 
\begin{align*} f(x_1, \cdots, x_v) &= \sum_{j=1}^{v} \alpha_j x_j + G(x_1, \cdots, x_v) \\ 
\text{where} |G(x)| &\leq C \sum_{j=2}^{v}(x_j - x_1)^2.
\end{align*} 
\end{mythm*}
\end{quote}
The local field version of this theorem is almost exactly the same:
\begin{mythm}\label{simul-thm} 
Given any constant $C > 0$ and a vector $\alpha \in K^v$ such that 
\begin{equation}\label{alpha-hypothesis}
\alpha \cdot u \neq 0  \text{ for every } u \in \{ 0, 1\}^v \text{ with } u \ne 0,  u \neq (1, 1, \ldots, 1)
\end{equation}  
and such that
\begin{equation} \label{alpha-cancellation}
\sum_{j=1}^v \alpha_j = 0,
\end{equation}  
there exists a positive constant $c(\alpha, K)$ and a set $E = E(C, \alpha, K) \subseteq R$ of Hausdorff dimension $c(K, \alpha) > 0$ with the following property. 
\vskip0.1in
\noindent The set $E$ does not contain any nontrivial solution of the equation \[f(x_1, \cdots, x_v) = 0, \qquad (x_1, \cdots, x_v) \text{ not all identical}, \] for any twice strictly differentiable function $f$ of the form \begin{align} f(x_1, \cdots, x_v) &= \sum_{j=1}^{v} \alpha_j x_j + G(x_1, \cdots, x_v)  \label{f-form} \\ \text{  where } |G(x)| &\leq C \sum_{j=2}^{v}|x_j - x_1|^2.  \label{G}\end{align} 

Here, the absolute value on $K$ is chosen so that a uniformizing element $t$ satisfies $|t| = q$, where $q$ is the number of elements of the residue class field of the non-archimedean local field $K$.
\end{mythm} 
\section{The arithmetic of complete, nonarchimedean local fields}
\subsection{A height function for $\mathbb{Q}_p$ and for complete local fields of finite characteristic}
Let $K$ be a function field; that is, a local field of characteristic zero. Let $x \in R$, the ring of integers of $K$. We can write $x$ in the form
\[\sum_{j=0}^{\infty} x_j t^j\]
a polynomial in the variable $t$, where $x_j \in \mathbb{F}_{q}$ for some $q = p^f$. Addition in $R$ and multiplication in $R$ work in the usual way: The sum $x + y$ is 
\[\sum_{j=0}^{\infty} (x_j + y_j) t^j\]
and the product of $x$ and $y$ is
\[\sum_{j=0}^{\infty} \left(\sum_{k_1 + k_2 = j} x_{k_1} y_{k_2} \right)t^j.\]
Each of the sums $\sum_{k_1 + k_2 = j} x_{k_1} y_{k_2}$ is finite and is therefore well-defined in $\mathbb{F}_q$.

The point is that, if we know that $x$ and $y$ are not just power series, but polynomials, then the degree of $x + y$ is bounded above by the maximum of the degrees of $x$ and $y$, and the degree of $xy$ is bounded above by the sum of the degrees of $x$ and $y$. We define a height function on $R$ in the following way: if $0 \neq x \in R$ is a polynomial then the height of $x$ will be defined to be the degree of $x$, and if $x$ is not a polynomial then the height will be defined to be $\infty$. 

Similarly, each element $x$ of $\mathbb{Z}_p$ can be written
\[x = \sum_{j=0}^{\infty} x_j p^j.\]
If the sum is finite, then $x \in \mathbb{Z}$ is an integer written in base $p$. So we can define a height on $\mathbb{Z}_p$ as follows: the height of $x \in \mathbb{Z}_p$ is the number of digits in the base $p$ expansion of $x$ minus one if $x$ has a finite expansion in $\mathbb{Z}_p$, and infinity otherwise.

Unfortunately, this height function does not behave quite as well as the one on $\mathbb{F}_{q}[[x]]$: the sum of two numbers of height $h_1, h_2$ may have height as large as $\max(h_1, h_2) + 1$. The product, however, is well-behaved: if $x$ and $y$ have height $h_1$ and $h_2$ then the product $xy$ has height at most $h_1 h_2$.

Given any complete, nonarchimedean local field $K$, we seek to define a height function on $R$, the ring of integers of $K$, satisfying the following properties:
\begin{itemize}
\item There are only finitely many elements of $R$ with height at most $h$ for any finite number $h$.
\item If $x$ has height $h_1$ and $y$ has height $h_2$, then $x + y$ has height at most $\max(h_1, h_2) + C$ for some $C$ depending only on the local field $K$.
\item If $x$ has height $x_1$ and $y$ has height $h_2$, then $xy$ has height at most $h_1 + h_2 + C$ for some $C$ depending only on the local field $K$.
\end{itemize}
We have already constructed such a function for every complete, nonarchimedean local field of finite characteristic. We therefore only need to construct a height function for finite extensions of $\mathbb{Q}_p$.
\subsection{Defining a height function for unramified extensions of $\mathbb{Q}_p$}
We first consider an unramified extension $K/\mathbb{Q}_p$. $K$ is formed by enlarging the residue field of $\mathbb{Q}_p$. Let $\mathbb{F}_{p^f}$ be the residue field of $K$. That is, the field $R/pR$ is isomorphic to $\mathbb{F}_{p^f}$. Then, select $\alpha \in R$ such that $\alpha \text{ (mod $pR$)}$ generates the multiplicative group $\mathbb{F}_{p^f}$. Then $1, \alpha, \alpha^2, \ldots, \alpha^{f-1}$ are linearly independent over $\mathbb{Q}_p$. Notice that $\alpha$ satisfies the relation $b(\alpha) \equiv 0 \text{ (mod $pR$)}$ where $b(\alpha)$ is the $p^f$th cyclotomic polynomial. Furthermore, the derivative $b^{\prime}(\alpha)$ is seen to have absolute value $1$, because the cyclotomic polynomial $b$ on $\mathbb{F}_{p^f}$ does not have any multiple roots. Therefore, we can apply Hensel's lemma to conclude that there is an element $t \in R$ such that $b(t) = 0$. We know that $\{1, t, t^2, \ldots, t^{f-1}\}$ is a basis for $\mathbb{Q}_p$ since reducing mod $p$ gives a basis for $\mathbb{F}_{p^f}$. Notice that each coefficient of $b$ is an integer, and therefore can be viewed as an element of $\mathbb{Z}_p$ of finite height. Therefore, each power $\{t^j : 0 \leq j < 2f -1\}$ can be written as a $\mathbb{Z}_p$-linear combination of $1, t, t^2, \ldots, t^{f-1}$ where each element has finite height.

Because $\{1, t, t^2, \ldots, t^{f-1}\}$ form a basis of the free module $R/\mathbb{Z}_p$, we can write each element $x \in R$ in the form
\[x^{(0)} + t x^{(1)} + \cdots + t^{f-1} x^{(f-1)}\]
where $x^{(k)} \in \mathbb{Z}_p$ for all $j,k$. We define the height of an element of $R$ to be the maximum of the heights of $x^{(0)}, \ldots, x^{(f-1)}$, if it is finite. It immediately follows that if $x$ and $y$ have heights $h_1$ and $h_2$, then $x + y$ has height no more than $\max(h_1, h_2) + 1$ by the behaviour of the height function on $\mathbb{Q}_p$.

As for the product, consider $xy$ where $x$ has height $h_1$ and $y$ has height $h_2$. Then $xy$ is
\[xy = \sum_{j = 0}^{f-1} \sum_{k=0}^{f-1} x_j y_k t^{j+k}.\]
We observed above that $t^{j+k}$ can be written as a linear combination of $1, \ldots, t^{f-1}$ of finite height, so $xy$ can be written as a sum a sum of a bounded (that is, depending only on the field $K$) number of terms that have height at most $h_1 + h_2 + C^{\prime}$, where $C^{\prime}$ is a constant depending only on the field $K$. Therefore, the product $xy$ has height at most $h_1 + h_2 + C$, where $C$ is a constant depending only on the field $K$.
\subsection{Defining a Height Function for Arbitrary Finite Extensions of $\mathbb{Q}_p$}
Now, let $K/\mathbb{Q}_p$ be an arbitrary finite extension of $\mathbb{Q}_p$. There is an intermediate field $L$, a maximal unramified subextension of $K/\mathbb{Q}_p$, such that $L/\mathbb{Q}_p$ is an unramified extension and $K/L$ is a totally ramified extension.

We use the following basic fact about totally ramified extensions of $L$: every totally ramified extension $K$ of $L$ is generated by a root $s$ of an Eisenstein polynomial over $R_L$, the ring of integers of $L$. A proof of this fact can be found in \cite{R00}, Chapter 2, Section 4.2 for the special case $L = \mathbb{Q}_p$, but the proof extends to arbitrary $L$. This root $s$ can be chosen to satisfy $|s| = q^{-1}$, i.e. $s$ can be taken to be a uniformizer for $K$. Let $a(x) = x^e + a_{e-1} x^{e-1} + \cdots + a_1 x + a_0$ be this Eisenstein polynomial. Then each of $a_{e-1}, \ldots, a_0 \in R_L$ is divisible by $p$, with $|a_0|_{L} = q^{-1}$. Thus $|a_0|_K = q^{-e}$ because we normalized the absolute value so that $|s| = q^{-1}$.

Consider the equation $a(s) = 0$, which holds in $K$. Expanding the left side of the equation, we get $s^e + a_{e-1}s^{e-1} + \cdots + a_1 s + a_0 = 0$. The derivative $a^{\prime}$ is nonzero at $s$; otherwise, $s$ would have degree less than $e$ over $K$, which is impossible because the polynomial $a$ is irreducible by Eisenstein's criterion. Suppose that the absolute value of $a^{\prime}(s)$ in $K$ is equal to $q^{-\alpha}$. Reduce the equation $a(t) = 0$ mod $s^{2 \alpha + 1} R$. Performing this reduction we get $s^e + a_{e-1} s^{e-1} + \cdots + a_1 t + a_0 \equiv 0$ (mod $s^{2 \alpha + 1}R)$. Evidently this equation continues to hold if we replace the constants $a_0, \ldots, a_{e-1}$ by any other constants that are congruent to $a_0, \ldots, a_{e-1}$ mod $s^{2 \alpha + 1}R$; in particular, we can replace them with elements $b_0, \ldots, b_{e-1}$ of $L$ with $L$-height no more than $\lceil \frac{2 \alpha + 1}{e} \rceil$. Let $b(x)$ be a polynomial with this replacement made. Then $b(s) \equiv 0$ (mod $s^{2 \alpha + 1} R$) and $|b^{\prime}(t)| = q^{\alpha}$, so by the version of Hensel's lemma appearing in Chapter 2, section 1.5 of \cite{R00}, it follows that $b(x)$ has a root within a $p^{-\alpha}$-neighbourhood of $s$. In particular, this root must have absolute value $q^{-1}$ because $\alpha + 1 \geq 1$. Let $t$ be this root of $b(x)$. Because $t$ is a uniformizer of $K$, it follows that $1, t, t^2, \ldots, t^{e-1}$ form a basis for $K/L$, we can write every element of the ring of integers of $K$ in the form
\[x = x^{(0)} + x^{(1)}t + \cdots + x^{(e-1)}t^{e-1}\]
where each of $x^{(0)}, \ldots, x^{(e-1)}$ are in $L$. Define the height of $x$ to be the maximum of the heights of $x^{(0)}, \ldots, x^{(e-1)}$. Suppose $x, y \in K$ such that $x$ has height $h_1$ and $y$ has height $h_2$. Then $x + y$ has height at most $\max(h_1, h_2) + 1$, which follows from the bound on the heights of $x^{(0)} + y^{(0)}, \ldots, x^{(e-1)} + y^{(e-1)}$ obtained in the previous section. 

It remains to be seen that $xy$ has height $h_1 + h_2 + C$, where $C$ is a constant depending only on the field $K$. 

Once again, we have
\[xy = \sum_{j = 0}^{f-1} \sum_{k=0}^{f-1} x_j y_k t^{j+k}.\]
The number of summands depends only on the field $K$ (and not on $h_1$ or $h_2$) and the height of each term is bounded above by $h_1 + h_2 + C^{\prime}$, where $C^{\prime}$ is a constant depending only on the field $K$, because the coefficients of $b$ have finite height. Thus the height of $xy$ is bounded above by $h_1 + h_2 + C$, where $C$ is a constant depending only on the field $K$.

Combining these facts with the previous subsection, we arrive at a basis $\{t_1^{k_1} t_2^{k_2}\}$ of the extension $K/\mathbb{Q}_p$, where $k_1$ runs from $0$ to $f$ and $k_2$ runs from $0$ to $e$, with respect to which the height is defined. Here $|t_1| = 1$ and $|t_2| = q^{-1} := p^{-f}$, because $t_1$ has absolute value $1$ in $L$ (and therefore in $K$) and $t_2$ is a uniformizing element of $K$. We also have $|p| = |t_2|^e = q^{-e}$. 
\subsection{Negatives of Elements of $K$ with finite height}
In the case of the ring $\mathbb{F}_q[[t]]$, the negative of an element of height $h$ will also have height $h$. Unfortunately, for other discrete valuation rings, the negative of an element of finite height will typically have infinite height. However, we can still say something useful about such elements.

Every element $z$ of $\mathbb{Z}_p$ has an expansion
\[\sum_{j=0}^{\infty} z_j p^j.\]
We will call the collection of digits $\{z_0, \ldots, z_{d-1} \}$ the $d$ least significant digits of $z$. We will say that $z^{(1)}$ and $z^{(2)} \in \mathbb{Z}_p$ differ only in the $d$ least significant digits if $z^{(1)}_j = z^{(2)}_j$ for all $j \geq d$. We will extend this notion to other fields as follows.

Let $K$ be a complete, nonarchimedean local field of characteristic zero, and $R$ be the ring of integers of $K$. Then, using the basis $\{t_1^{k_1} t_2^{k_2}\}$ described in the previous subsection, we can write an arbitrary element of $R$ in the form
\[x = \sum_{k_1 = 0}^{f-1} \sum_{k_2 = 0}^{e-1} x^{(k_1, k_2)} t_1^{k_1} t_2^{k_2}\]
where $x^{(k_1, k_2)}$ lie in $\mathbb{Z}_p$. We will say that $x$ and $y$ differ only in the $d$ least significant digits if $x^{(k_1, k_2)}$ and $y^{(k_1, k_2)}$ differ only in the $d$ least significant digits for all $k$.

We now state the following lemma:
\begin{mylem}\label{neg-lem}
Let $K/\mathbb{Q}_p$ be a finite extension with inertia degree $f$ and ramification index $e$, and let $R$ be the ring of integers of $K$. Then there exists a set $\{a_1, \ldots, a_{2^{ef}}\}$ of elements of $R$ such that if $x \in R$ has height at most $h$, then $-k$ differs from one of $a_1, \ldots, a_{2^{ef}}$ in only at most the first $h$ digits.
\end{mylem}
\begin{proof}
For $x \in \mathbb{Z}_p$ of height at most $h$, it is easy to see that either
\[-x = \begin{cases}
\sum_{j=0}^{h-1} y_j p^j + \sum_{j=h}^{\infty} (p-1) p^j & \text {if $x \neq 0$} \\
0 & \text{if $x = 0$}
\end{cases}\]
for some appropriate digits $y_0, \ldots, y_{h-1}$. Thus $x$ differs from one of $-1$ and $0$ only in the first $h$ digits. Therefore, for $x$ in $R$, we have that $-x^{(k_1, k_2)}$ differs from either $0$ or $-1$ in only the first $h$ digits. Therefore, we have
\[-x = \sum_{k_1 = 0}^{f-1} \sum_{k_2 = 0}^{e-1} - x^{(k_1, k_2)}  t_1^{k_1} t_2^{k_2}\]
differs from one of the expressions
\[\sum_{k_1 = 0}^{f-1} \sum_{k_2 = 0}^{e-1} \epsilon_{k_1, k_2} \cdot (-1) t_1^{k_1} t_2^{k_2}\]
in at most the first $h$ digits, for $\epsilon_{k_1, k_2} = 0$ or $1$.
\end{proof}
The point of this lemma is to establish the following corollary:
\begin{mycor}\label{perturbation}
Let $x$ and $y$ be elements of $R$ with height at most $h$, and let $\delta \in R$, $\delta \neq 0$ satisfy the inequality $|\delta| \leq q^{-e (h+1)}$. Then we have that $|(-x + \delta) - y| \geq |\delta|$.
\end{mycor}
\begin{proof}
We first consider the case in which $x = 0$. If $y$ is also equal to zero, then $(-x + \delta) - y$ is simply $\delta$ and clearly $|\delta| \geq |\delta|.$

If $x = 0$ and $y$ is nonzero, then $|(-x + \delta) - y|$ is equal to $|y - \delta|$. Since $y$ is nonzero and has height at most $h$, it follows that $|y| \geq q^{-e(h-1) -(e - 1)} = q^{-eh + 1}$. Because $|\delta| \leq q^{-e(h+1)} < q^{-eh + 1}$ it follows from the ultrametric inequality that $|y - \delta| = |y| > |\delta|$.

So we are left with the case in which $x \neq 0$. Because $x$ has height at most $h$, it follows from Lemma \ref{neg-lem} that $-x$ differs from one of the nonzero elements $a_1, \ldots, a_{2^{ef}}$ in only the first $h$ digits. In particular, this implies that if we write
\[x = \sum_{0 \leq k_1 < f} \sum_{0 \leq k_2 < e} \left(\sum_{j=0}^{\infty} x_j^{(k_1, k_2)} p^j  \right) t_1^{k_1} t_2^{k_2}\]
then $x_h^{(k_1,k_2)}$ must be equal to $p-1$ for at least one $(k_1, k_2)$. In particular, $x_h^{(k_1, k_2)}$ is nonzero. Thus, the same is true for $-x + \delta$ because $|\delta| \leq q^{-e(h+1)}$, and each term $(-x)_j^{(k_1, k_2)}p^ht_1^{k_1} t_2^{k_2}$ such that $(-x)_h^{(k_1, k_2)}$ is nonzero has absolute value exactly $q^{-eh - k_2}$. However, $y_h^{(k_1, k_2)}$ is certainly zero because $y$ has height at most $h$. Therefore, $(- x + \delta)_h^{(k_1, k_2)}$ differs from $y_h^{(k_1, k_2)}$ for some $(k_1, k_2)$ and thus $|(-x + \delta) - y| \geq q^{-eh - k_2} > q^{-e(h+1)} \geq |\delta|$.
\end{proof}
\subsection{A set avoiding the zeros of a polynomial}
\begin{mypro}\label{mathe-mainprop}
Let $T_1, \ldots, T_{v-1}, T_v$ be sets, each of which is a union of balls of radius $\frac{1}{q^{\mu}}$, and let $p(x_1, \ldots, x_v)$ be a polynomial satisfying $\left|\frac{\partial p}{\partial x_v}\right| \geq q^{-A}$ for $(x_1, \ldots, x_v) \in T_1 \times \cdots \times T_v$. Let $D = n/d$ be the dimension given in the statement of Theorem \ref{mathe}. Then there exists a positive real number $c$ depending on $A$, the field $K$, and $\mu$ with the following property: Subdivide each such ball into balls of radius $\frac{1}{q^{\nu}}$, where $\nu > \mu$ is divisible by the ramification index $e$ of the extension $K/\mathbb{Q}_p$ if $K$ has characteristic zero. Then, if $\mu$ is sufficiently large, there exist sets $S_1 \subset T_1, \ldots, S_v \subset T_v$ such that:
\begin{enumerate}[(a)]
\item There are no solutions to $p(x_1, \ldots, x_v) = 0$ with $x_1 \in S_1, \ldots, x_v \in S_v$. Furthermore, $p$ satisfies the bound  $|p(x_1, \ldots, x_v)| \geq c q^{-\nu}$ on $S_1 \times \cdots \times S_v$.
\item For $i = 1, \ldots, v-1$ and any ball $U$ of radius $\frac{1}{q^{\nu}}$ contained in one of the $\frac{1}{q^{\mu}}$-balls of $T_i$, we have that $S_i \cap U$ is a ball of radius $\geq c q^{-\nu n/D}$.
\end{enumerate}
\end{mypro}
\begin{proof}
Let $h := \nu/e$, where $e = 1$ if $K$ has finite characteristic and $e$ is the ramification index of $K/\mathbb{Q}_p$ if $K$ has characteristic zero. It follows from our assumptions that $h$ is an integer. Notice that each ball $B$ of radius $q^{-\nu}$ contains exactly one element of $R^n$ of height at most $h$ (the height of an element of $R^n$ will be defined to be the maximum of the heights of the components.)

Suppose $x_1, \ldots, x_v \in R^n$, and each of $x_1, \ldots, x_v$ has height at most $h$. Let $p$ be an $nv$-variate polynomial of degree at most $d$. Let $s := \binom{d+nv}{nv}$ be the dimension of the space of $nv$-variate polynomials of degree $d$. Then $p(x_1, \ldots, x_v)$ is a sum of at most $s$ terms, each of which is either of height at most $b + hd + C$, or the negative of an element of height at most $b + hd + C$, where $C$ is a value that depends on the field $K$ and degree $d$ of $p$, but not on the height $h$. Thus we can write $p(x_1, \ldots, x_v)$ as a sum of $2$ terms: a term with height at most $b + dh + C + s$ and the negative of such a term. It immediately follows from Corollary \ref{perturbation} that if $|\delta| \leq q^{-e(b + dh + C + s + 1)}$ that $|p(x_1, \ldots, x_v) + \delta| \geq |\delta|$. 

This fact will serve as a local-field substitute for the following algebraic fact: if $p$ is a $nv$-variate polynomial of degree $d$ with coefficients in $\mathbb{Z}$, and $x_1, \ldots, x_v$ are multiples of $\frac{1}{N}$, then $p(x_1, \ldots, x_v)$ is a multiple of $\frac{1}{N^d}$, so $p(x_1, \ldots, x_v) + \delta$ differs from zero by at least $|\delta|$ so long as $|\delta| < \frac{1}{2N^d}$. This simple algebraic fact is the key to M\'ath\'e's construction in \cite{M17}.

The only remaining piece of the puzzle is the familiar fact that if a polynomial $p$ satisfies $q^{-A} \leq \left|\frac{\partial p}{ \partial x_v}\right|$ uniformly on a compact set $T_1 \times \cdots \times T_v$, then $|q^{-A} \delta| \leq |p(x_1, \ldots, x_v + \delta) - p(x_1, \ldots, x_v)| \leq |\delta|$ for all $\delta$ sufficiently small (recalling here that any polynomial with coefficients in $R$ has derivatives bounded in absolute value by $1$ on $R$). This implies, in light of the above arguments, that if $|\delta|$ is smaller than $q^{-A - e(b + dh + C + s + 1)}$, and $x_1, \ldots, x_v$ all have height at most $h$ where $h$ is sufficiently large depending on the polynomial $p$, then $p(x_1, \ldots, x_v + \delta)$ is not within $q^{-A} |\delta|$ of $0$.

Finally, since the partial derivatives of $p$ are all bounded above by $1$ in absolute value, we have that the same bound holds provided that $x_1, \ldots, x_v$ are all within closed balls of radius $q^{-A - e(b + dh + C + s + 1) - 1}$ determined in the following way: given $B \subset T_j$ for $1 \leq j \leq v-1$, we select $S_j \cap B$ to be the unique all containing an element of height $h$ Selecting $\delta$ with $|\delta| = q^{-A - e(b + dh + C + s + 1) - 1}$, we define $S_v \cap B$ To be the unique ball containing $x_B + \delta$, where $x_B$ is the unique element of height $h$ in $T_v$.
\end{proof}
\section{Avoidance of Configurations at a Single Scale}
Proposition \ref{mathe-mainprop} will be central to the construction described in Theorem \ref{mathe}.  In order to prove Theorems \ref{non-simul-thm} and \ref{non-simul-thm vector-valued}, we will need $p$-adic analogues of the lemmas appearing in \cite{FP16}. The first step in proving such an analogue is counting the balls intersected by the zero set of $f$. Notice that for $m = n = 1$, the condition that $Df$ is of full rank is equivalent to the statement that $f$ has a nonvanishing partial derivative.
\begin{mylem}\label{total-box-lemma}
There exists $\mu_0(C_0, C_1, C_2) > 0$ such that the following statement holds for all $\mu \geq \mu_0$. Let $\mathbb{T}$ be an $nv$-dimensional ball of radius $q^{- \mu}$ and let  $f(x_1, \ldots, x_v) : \mathbb{T} \to K^m$ be a function such that $Df$ has an $m$-by-$m$ minor that is, in absolute value, at least a constant $C_0$ on all of $\mathbb{T}$, and whose entries are bounded above in absolute value by a constant $C_1$. Suppose further that $C_2$ is an upper bound for the operator norm of the second derivative of $f$. Let $Z_f$ be the set of $(x_1, \ldots, x_v)$ such that $f(x_1, \ldots, x_v) = 0$. Subdivide the ball $\mathbb{T}$ into balls of radius $q^{- \lambda}$. If $\lambda$ is sufficiently large, then the number of balls that intersect the zero set of $f$ is at most $C_3 q^{ - \mu + \lambda(nv - m)}$, where $\mu_0$ and $C_3$ depend on $C_0, C_1,$ and $C_2$ but not on $\lambda$.
\end{mylem}
\begin{proof}
Let $x_0$ be a point in $Z_f$. If $\mu$ is sufficiently large, then there exists a fixed $m$-by-$m$ submatrix $B$ of $Df$, indexed by the columns $(j_1, \ldots, j_m)$, such that $B$ has determinant with absolute value at least $C_0$ for all $x \in \mathbb{T}$. Consider the vector space $U$ spanned by the vectors $e_{j_i}$ whose $j_i$ component is $1$ and whose other components are $0$. Consider points of the form $x_0 + u$ where $u$ is a vector in the vector space $U$ of magnitude at most $q^{- \mu}$. By the assumptions on $Df$, we have $f(x_0 + u) = f(x_0) + Df_{x_0} u + O(\norm{u}^2)$. Here, $f(x_0) = 0$, and $\norm{Df_{x_0} u} \geq k_0 \norm{u}$, where $k_0 = C_0 C_1^{-(m-1)}$, as can be seen from the adjugate formula for the inverse of $B$. So we have that $\norm{f(x_0 + u)} \geq C_0 \norm{u} - O(\norm{u}^2)$. This is guaranteed to be positive provided that $\mu$ is sufficiently large.

Let $\lambda$ be larger than $\mu$. For a given $x \in Z_f$, consider the slab consisting of points $x + U + w$ where $w$ has a zero in each of the $j_i$ components and satisfies $|w| \leq q^{- \lambda}$; i.e., the $q^{- \lambda}$-neighbourhood of $x + U$. Let $x + u + w$ be some point in this slab. We have that $f(x + u + w) = f(x) + Df_{x}(u + w)  + O(\norm{u + w}^2) = Df_{x}(u + w) + O(\norm{u + w}^2)$. This has norm at least $(k_0 - C_2 q^{- \mu}) \norm{u} - C_1 \norm{w} + O(\norm{u + w}^2)$. Therefore, $f(x+u+w)$ is nonzero provided that $\norm{u}$ is at least $\frac{C_1}{k_0 - C_2 q^{- \mu} } \norm{w}$, and provided that $\mu$ is sufficiently large that the $O(\norm{u + w}^2)$ term has norm smaller than $\norm{u}$. Recall that $\norm{w} \leq q^{- \lambda}$, so if $\mu$ is sufficiently large this inequality will hold provided that $\norm{u}$ is at least $\frac{2C_1}{k_0} q^{- \lambda}$. Thus, subdividing this slab into $q^{-\lambda}$-balls, we have $Z_f$ will intersect only at most $\left(\frac{2 C_1}{k_0} \right)^{m}$ balls in the slab. Taking the union over disjoint parallel slabs that cover the $q^{- \mu}$ ball, we have a total of $O \left(\left(\frac{C_1}{k_0}\right)^{m} q^{\lambda(nv-m) - \mu}\right)$ intersections. 
\end{proof}
\begin{myrmk}
In the $m = n = 1$ case, we will not have any information about any of the second derivatives of $f$. Nonetheless, it is a simple consequence of the local-field version of the implicit function theorem that the equivalent of Lemma \ref{total-box-lemma} holds in this case.
\end{myrmk}
We can apply this counting lemma to prove the following proposition:
\begin{mypro}\label{non-simul-mainprop}
Let $T_1, \ldots, T_v$ be disjoint sets that can be expressed as the union of balls of radius $q^{- \mu}$. Let $f(x_1, \ldots, x_v)$ be a function defined on $T_1 \times \cdots \times T_v$ with the property that the derivative $Df$ of $f$ has full rank, and has a minor bounded below in absolute value by $C_0$, entries whose absolute values are bounded above by $C_1$, and $f$ has a second derivative with operator norm bounded above by $C_2$ on $T_1 \times \cdots \times T_v$, and let $D = \frac{m}{n(v-1)}$ be the dimension given in the statement of Theorem \ref{non-simul-thm}. Then there exists a $\nu_0(\mu, C_0, C_1, C_2)$ such that for all $\nu > \nu_0$  there exist $S_1 \subset T_1, \ldots, S_v \subset T_v$ satisfying the following conditions:
\begin{enumerate}
\item There are no solutions to $f(x_1, \ldots, x_v) = 0$ with $x_1 \in S_1, \ldots, x_v \in S_v$.
\item For any $1 \leq j \leq v-1$, and each ball $B$ of radius $q^{- \nu}$ contained in $T_j$, we have that $B \cap S_j$ is a ball of radius at least $c q^{-\frac{n}{D} \nu - C \mu}$ for some constants $c(C_0, C_1, C_2)$ and a constant $C$ depending only on $n, m$, and $v$.
\item For each ball $B$ of radius $q^{- \nu}$ contained in $T_j$, except for at most a  $q^{- \mu}$-fraction, we have that $B \cap S_v$ is a union of balls of radius $c q^{-\frac{n}{D} \nu - C \mu}$ for some small constant $c(C_0, C_1, C_2)$ and a constant $C$ depending only on $n, m,$ and $v$.
\end{enumerate}
\end{mypro}
The proof of this proposition, given the lemma, is a purely combinatorial argument almost identical to the proof of Proposition 3.4 appearing in \cite{FP16}. Therefore, we will provide only a summary of the local field version of the argument.

The proof is based on a ``projection lemma," which is a small modification of Lemma 3.5 appearing in \cite{FP16}:
\begin{mylem}\label{proj-lemma}
Let  $\mu, \nu, \lambda$ be such that $\mu \ll \nu \ll \lambda,$ and $\lambda$ satisfies $q^{- \lambda} = c q^{-n/D \nu - C \mu}$ for some appropriate constants $c$ and $C$. Let $T \times T^{\prime}$ be a union of $nr$-dimensional balls of radius $q^{-\mu}$, where $T$ is $n$-dimensional and $T^{\prime}$ is $n(r-1)$ dimensional. The collection of balls that constitute $T$ will be denoted $\mathbb{T}$, and the collection of balls that constitute $T^{\prime}$ will be denoted $\mathbb{T}^{\prime}$. Let $\mathbb{B} \subset \mathbb{T} \times \mathbb{T}^{\prime}$ be a collection of balls of radius $q^{- \lambda}$, where $\lambda$ is sufficiently large depending on $\mu$, whose union is denoted $B$, and let $\nu$ be such that $\mu \ll \nu \ll \lambda$. Then there exist sets $S \subset T$ and $\mathbb{B}^{\prime} \subset \mathbb{T}^{\prime}$ such that
\vskip0.1in 
\begin{enumerate}[(a)]
\item \label{S-construction} Let $\mathbb{T}_{\nu}$ be the collection of balls of radius $q^{-\nu}$. Then there exists a subset $\mathbb{T}_{\nu}^*$ of $\mathbb{T}_{\nu}$ with cardinality at least $q^{(\nu - 2 \mu)n}$ such that, for each $U \in \mathbb{T}_{\nu}^*$, the intersection $U \cap S$ consists of exactly one ball of radius $q^{- \lambda}$. For balls $U \notin \mathbb{T}_{\nu}^*$, the intersection $U \cap S$ will be empty.
\vskip0.1in   
\item \label{cardinality-B'} $\#(\mathbb{B}^{\prime}) \leq q^{\mu(n+1) + \nu n  - \lambda n} \#(\mathbb B)$.
\vskip0.1in
\item \label{set-containment} $(S \times T^{\prime}) \cap B \subseteq S \times B^{\prime}$, where $B^{\prime}$ is the union of the balls in $\mathbb{B}^{\prime}$.
\end{enumerate}   
\end{mylem}
This lemma allows us to project the $q^{-\lambda}$-neighbourhood of $Z_f$ onto sets of successively smaller dimension, using the sets $S$ from the lemma as the sets $S_1, \ldots, S_{v-1}$ promised by Proposition \ref{non-simul-mainprop}. After $r-1$ applications of Lemma \ref{proj-lemma}, we arrive at a collection of $n$-dimensional balls whose complement satisfies the conditions on $S_v$ of Proposition \ref{non-simul-mainprop}. The details of the strategy, and of the proof of Lemma \ref{proj-lemma}, are omitted because of their similarity to those of Lemma 3.4, Proposition 3.5, and Proposition 3.6 in Section 3 of \cite{FP16}.
\section{The construction of the sets in Theorems \ref{mathe}, \ref{non-simul-thm}, and \ref{non-simul-thm vector-valued}} 
\subsection{General Procedure}
The proofs of Theorems \ref{mathe}, \ref{non-simul-thm}, and \ref{non-simul-thm vector-valued} from the propositions above proceed in almost the same way as in \cite{FP16}. We will keep a running queue $\mathcal{Q}$ that will guide the construction. In the case of Theorems \ref{mathe} and \ref{non-simul-thm}, we choose a sequence of ``privileged" differential operators $D_1, \ldots, D_{r_{\ell}}$, where $D_{r_{\ell}}$ represents a partial $r_{\ell}$th derivative of $f$ that is known to be nonvanishing on $B^{nv_{\ell}}$. Such a partial derivative exists for Theorem \ref{non-simul-thm} by assumption and for Theorem \ref{mathe} because a nonzero polynomial always has some partial derivative that is constant and nonzero, provided that we make the additional assumption that the degree is less than $\text{char }R$ in the case where $\text{char }R$ is finite. Let $E_0$ be the ball $B^n$. 

One minor difference in the proof of Theorem \ref{mathe} is that the number of variables $v$ is allowed to vary. For the cases of Theorems \ref{non-simul-thm} and \ref{non-simul-thm vector-valued}, take $v = v_1 = v_2 = \cdots$ throughout the rest of this argument.

Fix a sequence $\epsilon_j$ such that $\epsilon_j \to 0$. We select $\lambda_0$ sufficiently large so that the ball $E_0$ contains at least $v_1 + 1$ balls of radius $q^{- \lambda_0}$. Let $B_1^{(0)}, \ldots, B_{M_0}^{(0)}$ be an enumeration of the balls of radius $q^{- \lambda_0}$ contained in $B^n$ and let $\Sigma_0$ be the family of $v_1$-tuples of distinct such balls, ordered lexicographically and identified in the usual way with the family of injections from $\{1, \ldots, v_1\}$ into $\{1, \ldots, M_0\}$. Let $\mathcal{Q}_0$ be the queue consisting of $4$-tuples
\[\{(1, k, \sigma, 0) : 0 \leq k \leq r_1 - 1, \sigma \in \Sigma_0\},\]
where the queue elements are ordered so that $(1, k, \sigma, 0)$ precedes $(1, k^{\prime}, \sigma^{\prime}, 0)$ whenever $\sigma < \sigma^{\prime}$ and $(1, k, \sigma, 0)$ precedes $(1, k^{\prime}, \sigma, 0)$ whenever $k > k^{\prime}$.

\paragraph{Stage 1.}
At Stage $1$, we consider the first queue element $(1, k, \sigma, 0)$. Let $T_1 = B_{\sigma(1)}^{(0)}$ , \ldots, $T_v =  \ B_{\sigma(v)}^{(0)}$. 

First, we will consider the case of Theorems \ref{mathe} and \ref{non-simul-thm}. Let $f = D_k f_1$. By the ordering of $\mathcal{Q}_0$, we know that $k = r_1 - 1$, and therefore, $\frac{\partial f}{\partial x_i} = D_{k+1}f_1$ is nonzero for some $i$. By compactness, we therefore have that there is a lower bound $q^{-A_1}$ on the derivative $\left| D_{k+1} f_1 \right|$. We decompose each ball of radius $q^{- \lambda_0}$ into balls of radius $q^{-\mu_1}$, where $\mu_1 \geq \lambda_0$ and $\mu_1 > \mu^*$, where $\mu^*$ is as in Proposition \ref{mathe-mainprop} or \ref{non-simul-mainprop} applied to $T_1, \ldots, T_v$. We apply either Proposition \ref{mathe-mainprop} or Proposition \ref{non-simul-mainprop} to arrive at sets $S_1, \ldots, S_v$ with the properties guaranteed by the corresponding proposition. We can select $\nu = \nu_1$ in the proposition to be sufficiently large that the quantity $q^{- \lambda_j} := c q^{- \nu_1 n /D}$ appearing in Proposition \ref{mathe-mainprop} or the quantity $q^{- \lambda_1} := c q^{-\nu_1 n /D + C \mu_1}$ appearing in Proposition \ref{non-simul-mainprop} is larger than $q^{-\nu_1 \left(n/D + \epsilon_1\right) + n \mu_1}$.

Now, we will consider the case of Theorem \ref{simul-thm}. In the case of Theorem \ref{simul-thm}, we know that $D_{k+1}f_1 = Df_1$ has full rank on $T_1 \times \cdots \times T_v$ by assumption. Because this set is compact, it follows that there exists some $C_0$ such that, for each $x \in T_1 \times \cdots \times T_v$, $Df_1$ has a minor with absolute value greater than or equal to $C_0$ at $x$, and that each of the entries of $Df$ will be bounded above in absolute value by $C_1$ for some $C_1$. The uniform continuity of $D^2f$ also guarantees that $D^2f$ will be bounded above by $C_2$ in operator norm, for some appropriate value $C_2$. We can then select $\mu_1 \geq \mu^*$, where $\mu^*$ is as in Proposition \ref{non-simul-mainprop}. Select $\nu = \nu_1$ in Proposition \ref{non-simul-mainprop} to be sufficiently large so that that the quantity $q^{-\lambda_1} := c q^{- \nu_1/D + C \mu_1}$ appearing in Proposition \ref{non-simul-mainprop} is larger than $q^{-\nu_1 \left(n/D + \epsilon_1 \right) + n \mu_1}$.

In any case, we arrive at sets $S_1 \subset T_1, \ldots, S_v \subset T_v$ with the property that $D_kf_1$ is nonzero for $x_1 \in S_1, \ldots, x_v \in S_v$. We will define a subset $E_1 \subset E_0$ in the following way. We take  $E_1 \cap T_1 = E_0 \cap S_1$, $E_1 \cap T_2 = E_0 \cap S_2$, \ldots, $E_1 \cap T_v = E_0 \cap S_v$. All $x \in E_0 \setminus (T_1 \cup \cdots \cup T_v)$ will be in $E_1$. This gives a subset $E_1 \subset E_0$ that can be expressed as a disjoint union of balls of radius $q^{-\lambda_1}$.

Let $\mathcal{E}_1$ be the collection of balls of radius $q^{-\lambda_1}$ whose disjoint union is $E_1$. Enumerate the balls of $\mathcal{E}_1$ as $B_1^{(1)}, \ldots, B_{M_1}^{(1)}$. For $\ell = 1, 2$, define $\Sigma_1^{(\ell)}$ to be the collection of $v_{\ell}$ tuples of distinct such balls, ordered lexicographically and identified in the usual way with the family of injections from $\{1, \ldots, v_{\ell}\}$ into $\{1, \ldots, M_1\}$. We then form the queue $\mathcal{Q}_1^{\prime}$ consisting of $4$-tuples of the form
\[\{(\ell, k, \sigma, 1) : 1 \leq \ell \leq 2; 0 \leq k \leq r_{\ell} - 1, \sigma \in \Sigma_1^{(\ell)}\},\]
arranged so that $(\ell, k, \sigma, 1)$ precedes $(\ell^{\prime}, k^{\prime}, \sigma^{\prime}, 1)$ if $\ell \leq \ell^{\prime}$, so that $(\ell, k, \sigma, 1)$ precedes $(\ell, k^{\prime}, \sigma^{\prime}, 1)$ if $\sigma < \sigma^{\prime}$, and so that $(\ell, k, \sigma, 1)$ precedes $(\ell, k^{\prime}, \sigma, 1)$ if $k > k^{\prime}$. Arrive at the queue $\mathcal{Q}_1$ by appending the queue $\mathcal{Q}_1^{\prime}$ to $\mathcal{Q}_0$.
\paragraph*{Stage j.}
We will now describe Stage $j$ of the construction for $j > 1$. We follow essentially the same procedure as in Stage 1. We begin with a decreasing family of sets $E_0, \ldots E_{j-1}$. Each $E_{j^{\prime}}$ is a union of balls of radius $q^{-\lambda_{j^{\prime}}}$, the collection of which is called $\mathcal{E}_{j^{\prime}}$. The family of $v_{\ell^{\prime}}$ tuples of distinct balls in $\mathcal{E}_{j^\prime}$ will be denoted $\Sigma_{j^{\prime}}^{(\ell^{\prime})}$. We have a queue $\mathcal{Q}_{j-1}$ consisting of $4$-tuples $(\ell^{\prime}, k^{\prime}, \sigma^{\prime}, j^{\prime})$, where we have $0 \leq j^{\prime} \leq j-1$, $1 \leq \ell \leq j^{\prime} + 1$, $0 \leq k^{\prime} \leq r_{\ell^{\prime}} - 1$, and $\sigma^{\prime} \in \Sigma_{j^{\prime}}^{(\ell^{\prime})}$. The set $E_{j-1}$ has the property that $D_{k^{\prime}} f_{\ell^{\prime}}(x_1, \ldots, x_v) \neq 0$ for $x_1 \in B_{\sigma^{\prime}(1)}^{(j^{\prime})} \cap E_{j-1}, \ldots,  x_v \in B_{\sigma^{\prime}(v_{\ell^{\prime}})}^{(j^{\prime})} \cap E_{j-1}$ for any $(\ell^{\prime}, k^{\prime}, \sigma^{\prime}, j^{\prime})$ in the first $j-1$ elements of the queue $\mathcal{Q}_{j-1}$. Consider the $j$th queue element $(\ell, k, \sigma, j_0)$, where $\ell \leq j_0 \ll j$. Let $T_1, \ldots, T_v$ be the sets $B_{\sigma(1)} \cap E_{j-1}, \ldots, B_{\sigma(v_{\ell})} \cap E_{j-1}$. We consider a variety of cases depending on whether $k = r_{\ell} - 1$ or $k < r_{\ell} - 1$, and on whether we are considering Theorem \ref{mathe}, Theorem \ref{non-simul-thm}, or Theorem \ref{non-simul-thm vector-valued}.

\subparagraph{Case 1: Theorem \ref{mathe} or \ref{non-simul-thm}, $k = r_{\ell} - 1$.} Let $f = D_{r_{\ell} - 1} f_{\ell}$. In this case, we have that $k + 1 = r_{\ell}$; therefore, it follows by assumption that $D_{r_{\ell}}f_{\ell}$ is nonzero on all of $T_1 \times \cdots \times T_v$. By the definition of the differential operators $D_{r_{\ell}}$ and $D_{r_{\ell} - 1}$, we have that some partial first derivative of $f$ is nonvanishing on all of $T_1 \times \cdots \times T_v$. Let $q^{-A}$ be the lower bound on this partial derivative. Select $\mu_j \geq \lambda_{j-1}$ to be larger than the quantity $\mu^*$ appearing in Proposition \ref{mathe-mainprop} or Proposition \ref{non-simul-mainprop}, as is appropriate. We can then apply Proposition \ref{mathe-mainprop} or Proposition \ref{non-simul-mainprop} to the sets $T_1, \ldots, T_v$, with the quantity $\nu = \nu_j$ taken to be sufficiently large that the quantity $q^{- \lambda_j} := c q^{- n \nu_j/D}$ appearing in Proposition \ref{mathe-mainprop} or the quantity $q^{- \lambda_j} :=  c q^{- n \nu_j/D - C \mu_j}$ appearing in Proposition \ref{non-simul-mainprop} is larger than $q^{- \nu_j \left(1/D + \epsilon_j \right) + n \mu_j}$. 

\subparagraph{Case 2: Theorem \ref{mathe} or \ref{non-simul-thm}, $k < r_{\ell} - 1$.} If $k < r_{\ell} - 1$, then, by the ordering of the elements in the queue $\mathcal{Q}_{j-1}$, we will have that the $j-1$st element of $\mathcal{Q}_{j-1}$ is $(\ell, k + 1, \sigma, j_0)$. Therefore, by the previous stage, we have that for $x \in T_1, \ldots, x_v \in T_v$ that $D_{k+1} f_{\ell}$ is nonzero. But this implies that there exists some $A$ such that $D_{k+1}f$ is at least $q^{-A}$ in absolute value on all of $T_1 \times \cdots \times T_v$.  Select $\mu_j \geq \lambda_{j-1}$ to be larger than the quantity $\mu^*$ appearing in Proposition \ref{mathe-mainprop} or Proposition \ref{non-simul-mainprop}, as is appropriate. We can then apply Proposition \ref{mathe-mainprop} or Proposition \ref{non-simul-mainprop} to the sets $T_1, \ldots, T_v$ with the quantity $\nu = \nu_j$ taken to be sufficiently large that the quantity $q^{- \lambda_j} := c q^{- n \nu_j/D}$ appearing in Proposition \ref{mathe-mainprop} or the quantity $q^{-\lambda_j} := c q^{- \nu_j/D - C \mu_j}$ appearing in Proposition \ref{non-simul-mainprop} is larger than $q^{- \nu_j \left(1/D + \epsilon_j \right) + n \mu_j}$. 

\subparagraph{Case 3: Theorem \ref{non-simul-thm vector-valued}.} In this case, $k$ will always be chosen to be equal to $1$. Therefore, we have that on $T_1 \times \cdots \times T_v$ that $Df$ has a nonvanishing minor and that $D^2f$ is continuous. By compactness, we conclude that $Df$ is bounded below in absolute value by some value $C_0$ and the entries of $Df$ are bounded above in absolute value by some number $C_1$ on $T_1 \times \cdots \times T_v$. We also have an upper bound $C_2$ on the operator norm of $D^2f$. We can then select $\mu_j \geq \lambda_{j-1}$ to be larger than the quantity $\mu^*$ appearing in Proposition \ref{non-simul-mainprop}. We can then apply Proposition \ref{non-simul-mainprop} to the sets $T_1, \ldots, T_v$ with the quantity $\nu_j$ taken to be sufficiently large that the quantity $q^{-\lambda_j} := c q^{- \nu_j/D - C \mu_j}$ appearing in Proposition \ref{non-simul-mainprop} is larger than $q^{- \nu_j \left(n/D + \epsilon_j \right) + n \mu_j}$. 

In any case, we arrive at sets $S_1 \subset T_1, \ldots, S_v \subset T_v$ such that $D_{k} f_{\ell}$ is nonzero for $(x_1, \ldots, x_v) \in S_1 \times \cdots \times S_v$. We will define a subset $E_j \subset E_{j-1}$ in the following way. We take  $E_j \cap T_1 = E_{j-1} \cap S_1$, $E_j \cap T_2 = E_{j-1} \cap S_2$, \ldots, $E_j \cap T_v = E_{j-1} \cap S_v$. All $x \in E_{j-1} \setminus (T_1 \cup \cdots \cup T_v)$ will be in $E_j$. This gives a subset $E_j \subset E_{j-1}$ that can be expressed as a disjoint union of balls of radius $q^{-\lambda_j}$. Call the collection of such balls $\mathcal{E}_j$, and let $B_1^{(j)}, \ldots, B_{M_j}^{(j)}$ be an enumeration of the balls in $\mathcal{E}_j$. For each $1 \leq \ell \leq j,$ we define $\Sigma_j^{(\ell)}$ to be the collection of $v_{\ell}$-tuples of distinct balls in $\mathcal{E}_j$. We equip $\Sigma_j^{(\ell)}$ with the lexicographic order and identify $\Sigma_j^{(\ell)}$ with the set of injections from $\{1, \ldots, v_{\ell}\}$ into $\mathcal{E}_j$. Consider the queue $\mathcal{Q}_j^{\prime}$ consisting of $4$ tuples $(\ell, k, \sigma, j)$, where $1 \leq \ell \leq j+1$, $0 \leq k \leq r_{\ell} - 1$, and $\sigma \in \Sigma_j^{(\ell)}$. We order the queue $\mathcal{Q}_j^{\prime}$ in the following way: $(\ell, k, \sigma, j)$ will precede $(\ell^{\prime}, k^{\prime}, \sigma^{\prime}, j)$ if $\ell < \ell^{\prime}$, $(\ell, k, \sigma, j)$ precedes $(\ell, k^{\prime}, \sigma^{\prime}, j)$ if $\sigma < \sigma^{\prime}$, and $(\ell, k, \sigma, j)$ precedes $(\ell, k^{\prime}, \sigma, j)$ if $k > k^{\prime}$. We append the queue $\mathcal{Q}_j^{\prime}$ to $\mathcal{Q}_{j-1}$ to arrive at the queue $\mathcal{Q}_j$.
\subsection{Computation of Hausdorff Dimension}
We compute the Hausdorff dimension of the set $E$. Unlike the calculation in \cite{FP16}, we use the definition of Hausdorff dimension directly instead of appealing to Frostman's lemma.

Let $\mathcal{U}$ be a disjoint covering of a set $E$. Define the \emph{$s$-contribution} of a ball $V$ (not necessarily in $\mathcal{U}$) to be 
\[s(V) := \sum_{\substack{U \in \mathcal{U} \\ U \subset V}} r(U)^s\]
where $r(U)$ is the radius. Note that if $U \in \mathcal{U}$, then we have that $s(U) = r(U)^s$, and that if $V_1, \ldots, V_r$ are disjoint subsets of $V$, then $s(V_1) + \ldots + s(V_r) \leq s(V)$. 
\begin{mylem}\label{HausdorffLem}
Let $V$ be a ball of radius $q^{-\mu_k}$, and let $s < D$, where $D$ is the dimension promised for $E$ in Theorem \ref{mathe}, \ref{non-simul-thm}, or \ref{non-simul-thm vector-valued} Then:
\begin{enumerate}
\item If the majority of the volume of $\mathcal{U}$ contained in $V$ is in balls of $\mathcal{U}$ of radius strictly larger than $q^{- \mu_{k+1}},$ then $s(V) \geq \frac{1}{4} q^{- \mu_k s}$.
\item If the majority of the volume of $\mathcal{U}$ contained in $V$ is in balls of $\mathcal{U}$ of radius no more than $q^{- \mu_{k+1}}$, then $V$ contains t least $q^{(\mu_{k-1} - \mu_k)s}$ balls of radius $q^{- \mu_{k+1}}$ that contain a ball of $\mathcal{U}$.
\end{enumerate}
\end{mylem}
\begin{proof}[Proof of 1]
We consider two cases: the case in which
\begin{equation}\label{largeness}
\sum_{\substack{U \in \mathcal{U} \\ r(U) \geq q^{- \nu_j}}} r(U)^n \geq \frac{q^{- \mu_j n}}{4}
\end{equation}
and the case in which this inequality is not satisfied.
\paragraph*{Case 1}
If the inequality (\ref{largeness}) is satisfied, then we have
\begin{IEEEeqnarray*}{rCl}
\frac{q^{-\mu_j n}}{4} & \leq & \sum_{\substack{U \in \mathcal{U} \\ r(U) \geq q^{- \nu_j}}} r(U)^n\\
& = & \sum_{\mu_j \leq \rho \leq \nu_j} \#(\rho) q^{- n \rho}\\
& = & \sum_{\mu_j \leq \rho \leq \nu_j} \#(\rho) q^{-s \rho} q^{(s - n) \rho}\\
& \leq & \sum_{\mu_j \leq \rho \leq \nu_j} \#(\rho) q^{-s \rho} q^{(s - n)\mu_j}
\end{IEEEeqnarray*}
where $\#(\rho)$ is the number of balls of radius $q^{-\rho}$ in $\mathcal{U}$, and the last line holds because $s-n$ is negative. Dividing both sides of the inequality by $q^{(s - n) \mu_j}$ gives
\[\sum_{\mu_j \leq \rho \leq \nu_j} \#(\rho) q^{-s \rho} \geq \frac{q^{-s \mu_j}}{4}\]
as desired.
\paragraph*{Case 2} Suppose that the inequality (\ref{largeness}) is not satisfied. Let $\mathcal{U}^{\prime}$ be the collection of balls $U^{\prime}$ of radius $q^{-\nu_j}$ that contain a ball of $\mathcal{U}$. By Proposition 2.2 or Proposition 3.2, whichever is appropriate, we have that at least $\frac{1}{2}$ of the balls of radius $q^{-\nu_j}$ contained in $B$ intersect $E$, and since the balls of radius larger than $q^{- \nu_j}$ cover a set of measure strictly smaller than $\frac{q^{-\mu_j n}}{4}$, it follows that $\mathcal{U}^{\prime}$ must cover a set of measure at least $\frac{q^{- \mu_j n}}{4}$. Therefore $\mathcal{U}^{\prime}$ must consist of at least $\frac{q^{(\nu_j - \mu_j)n}}{4}$ balls.

Each ball in $\mathcal{U}^{\prime}$ intersects $E$ and therefore contains at least one ball in $\mathcal{E}_{j+1}$. Thus, letting $\mathcal{U}^{\prime \prime} \subset \mathcal{U}$ be the collection of balls in $\mathcal{U}$ that are contained in a ball in $\mathcal{U}^{\prime}$, we have that
\begin{IEEEeqnarray*}{rCl}
\sum_{U \in \mathcal{U}^{\prime \prime}} r(U)^s & \geq &  \frac{1}{4} q^{- (\frac{n}{D} \nu_j + C \mu_j) s} q^{\nu_j n - \mu_j n}\\
& = & \frac{1}{4} q^{\nu_j(n - \frac{n}{D}s) - \mu_j(n + Cs)}.\\
\end{IEEEeqnarray*}
Since we assumed $s < D$, we know that the coefficient on $\nu_j$ is strictly positive. Because $\nu_j > \exp(\mu_j)$ it follows that, for sufficiently large $j$, this is larger than $\frac{1}{4} q^{- \mu_j s}$. 
\end{proof}
\begin{proof}[Proof of 2]
This argument is similar to Case 2 above. We define $\mathcal{U}^{\prime}$ to be the collection of balls of radius $q^{- \mu_{j+1}}$ that contain a ball of $\mathcal{U}$. By the statements in Propositions 2.2 and 3.2, along with the fact that each $q^{- \lambda_{j+1}}$-ball that intersects $E$ has the property that each of the $q^{- \mu_{j+1}}$-balls that it contains intersects $E$,  we have that $V$ contains at least $q^{\nu_j - \mu_j}$ balls of radius $q^{- \mu_{j+1}}$ that intersect $E$, where $\mu_{j+1} = \frac{n}{D} \nu_j + C \mu_j$. Therefore, if at least half of the volume of $\mathcal{U}^{\prime}$ is contained in balls of radius at most $q^{- \mu_{j+1}}$, this means that there must be at least $q^{(\nu_j - \mu_j)n}/2$ balls of radius $q^{- \mu_{j+1}}$ that contain a ball of $\mathcal{U}^{\prime}$. But because $s > D$, it follows from the fact that $\nu_j \geq \exp(\mu_j)$ that $q^{\nu_j - \mu_j}/2 > q^{(\mu_{j+1} - \mu_j)s}$ for sufficiently large $j$, as in Case 2 above.
\end{proof}
From here, we can compute the Hausdorff dimension of the set $E$:
\begin{mypro}\label{HausdorffPro}
Let $V$ be a ball of radius $q^{-\mu_j}$ in $\mathcal{E}_j$, and let $\mathcal{U}$ be a covering of $V \cap E$ by balls of radius larger than $q^{-\mu_k}$ contained in $V$, where $k > j$. Let $s < D$. If $j$ is sufficiently large depending on $s$, then for the covering $\mathcal{U}$ we have $s(V) \geq \frac{1}{4}q^{-\mu_j s}$.
\end{mypro}
\begin{proof}
We prove this statement by induction on $k-j$. If $k-j = 1$, this statement is implied by part 1 of Lemma \ref{HausdorffLem}. So it only remains to show the inductive step.

Suppose first that the majority of the volume of $\mathcal{U}$ is contained in balls of radius strictly larger than $q^{- \mu_{j+1}}$. Then we can apply part 1 of Lemma \ref{HausdorffLem} to conclude that $s(V) \geq \frac{1}{4} q^{- \mu_j s}$, and we're done. Therefore, we can assume that the majority of the volume of $\mathcal{U}$ is contained in balls of radius at most $q^{- \mu_{j+1}}$. Then by part 2 of Lemma \ref{HausdorffLem}, there exist balls $V_1, \ldots, V_r$ of radius $q^{-\mu_{j+1}}$, where $r > q^{(\mu_{j+1} - \mu_j)s}$, such that each ball $V_j$ contains an element of $\mathcal{U}$. By the additivity of $s$ we have $s(V_1) + \ldots + s(V_r) \leq s(V)$, and by the inductive assumption we have that $s(V_t) \geq q^{- \mu_{j+1} s}$ for each $t$. Thus $s(V) \geq q^{- \mu_j s}$ as desired.
\end{proof}
The argument for the Minkowski dimension is slightly different.
\begin{mylem}
The Minkowski dimension of the set $E$ is equal to $1$ (in the case of Theorem \ref{non-simul-thm}) or $n$ (in the case of Theorems \ref{mathe} and \ref{non-simul-thm vector-valued}).
\end{mylem}
\begin{proof}
Let $N_{\mu}(E)$ be the number of closed balls of radius $q^{-\mu}$ that are required to cover $E$. Notice that this is equivalent to the number of balls of radius $q^{- \mu}$ that intersect $E$. So we need to count the number of such balls that intersect $E$. Let $q^{- \mu_j}$ denote the radius of the balls in $\mathcal{E}_j$.

If $\mu = \mu_j$ for some $j$, and $B$ is a ball of radius $q^{- \mu_{j-1}}$ that is not contained in $T_1, \ldots, T_v$ at stage $j$ of the construction, then the number of balls of radius $q^{- \mu_j}$ required to cover $B$ is precisely $\left(\frac{q^{- \mu_{j-1}}}{q^{- \mu_j}}\right)^n$. This is greater than or equal to $ q^{\mu_j^{-n + \epsilon}}$ (provided that $q^{\mu}$ is sufficiently small) by the growth rate of the $\nu_j$. 

Now suppose $\mu_j < \mu < \mu_{j+1}$, and consider any ball $B$ of radius $q^{-\mu_{j-1}}$ such that $B$ is not contained in $T_1 \cup \cdots \cup T_v$ at either stage $j$ or stage $j+1$ of the construction. Such a ball will always exist provided that $j$ is large enough.

Then by the argument from before, $E$ intersects all of the balls of radius $q^{- \mu}$ contained in $B$, and thus at least $q^{- \mu_j^{-n + \epsilon}}$ of radius $q^{- \mu_j}$ inside the ball $B$. Evidently $E$ must then intersect at least $\left(q^{- \mu_j + \mu} \right)^n \cdot q^{- \mu_j(-n + \epsilon)} \geq q^{\mu(n - \epsilon)}$ of the balls of radius $q^{-\mu}$, as desired.
\end{proof}
This completes the proofs of Theorems \ref{mathe}, \ref{non-simul-thm}, and \ref{non-simul-thm vector-valued}.
\section{Simultaneous avoidance of configurations with the same linearization} 
The other result applies to a (possibly uncountable) family of functions $f(x_1, \ldots, x_v)$ that have agreeing, unchanging, nondegenerate linearizations along the diagonal. The specific conditions are outlined in the statement of Theorem \ref{simul-thm}. 

As in the case of the Theorems \ref{non-simul-thm}, and \ref{non-simul-thm vector-valued}, the construction is a Cantor-like set constructed in the same manner as \cite{FP16}. This lemma outlines the strategy for a single step of the construction.

\begin{mylem}\label{simul-lem}
Let $B_1, B_2 \subset B$ be distinct (and thus disjoint) balls of radius $q^{- \lambda}$.  Let $\mathcal{A} \subset \{1, \ldots, v\}$ be an arbitrary nonempty proper subset of $\{1, \ldots, v\}$. Then we can find $B_1^{\prime} \subset B_1$, $B_2^{\prime} \subset B_2$ of radius at least $C_1 q^{- \lambda}$ such that $\abs{\alpha_1 x_1 + \cdots + \alpha_v x_v} \geq C_1r$ for $x_j \in B_1^{\prime}$ if $j \in \mathcal{A}$ and $x_j \in B_2^{\prime}$ if $j \notin \mathcal{A}$, where $C_1$ is a constant depending on the vector $\alpha$, the local field $K$, and on $\mathcal{A}$.
\end{mylem}
\begin{proof}
Since the bounds in this lemma are allowed to depend on $\alpha$, we will normalize $\alpha$ so that every component of $\alpha$ is in $R$ and at least one component of $\alpha$ is invertible in $R$; that is, the norm of $\alpha$ will be taken to be $1$. 

Assuming this normalization, let $C_1 = q^{-1}\abs{\sum_{j \notin \mathcal{A}} \alpha_j}$. This was assumed to be strictly positive in the statement of Theorem \ref{simul-thm}.

Then, selecting $B_1^{\prime}$ and $B_2^{\prime \prime}$ to be any balls of radius $C_1 q^{- \lambda}$ contained in $B_1$ and $B_2$, we consider the set $\alpha \cdot \mathcal{B}$, where $\mathcal{B} = B^{(1)} \times \cdots \times B^{(v)}$, and $B^{(j)} = B_1^{\prime}$ if $j \in \mathcal{A}$ and $B^{(j)} = B_2^{\prime \prime}$ otherwise. Since each component of $\alpha$ is in $R$, we have immediately that $\alpha \cdot \mathcal{B}$ is contained in a ball $A \subset R$ of radius at most $C_1 q^{- \lambda}$. If this ball is not the unique $C_1 q^{-\lambda}$-ball containing zero, then we can take $B_2^{\prime} = B_2^{\prime \prime}$ to complete the proof. If not, then select any element $b \in R$ with $|b| = q^{-\lambda - 1}$ and define $B_2^{\prime} := B_2^{\prime \prime} + b \subset B_2$. Then defining the ball $\mathcal{B}^{\prime}$ to be $B^{(1) \ast} \times \cdots \times B^{(v) \ast}$, with $B^{(j) \ast} = B_1^{\prime}$ if $j \in \mathcal{A}$ and $B^{(j) \ast} = B_2^{\prime}$ otherwise, we have that $\alpha \cdot \mathcal{B}^{\prime}$ maps into $A + b \sum_{j \in \mathcal{A}^C} \alpha_j$, which is a ball of radius $C_1 q^{-\lambda}$ that does not contain $0$.
\end{proof}
This lemma can be applied iteratively in the following way: Starting with a ball $B$ of radius $q^{-\lambda}$, we can pick two balls $B_1, B_2 \subset B$ of radius $q^{-\lambda - 1}$. Then, for some $\mathcal{A}$, apply Lemma \ref{simul-lem} to arrive $B_1^{\prime}$ and $B_2^{\prime}$. Now, pick a different $\mathcal{A}$ and apply the lemma to these balls to arrive at $B_1^{\prime \prime}$ and $B_2^{\prime \prime}$. Repeat this process for all nonempty proper subsets $\mathcal{A} \subset \{1, \ldots, v\}$, and we arrive at $B_1^{\ast} \subset B_1 , B_2^{\ast} \subset B_2$ where $B_1$ and $B_2$ have radius at least $C^* q^{-\lambda}$ for some $C^*$, and $\alpha_1 x_1 + \cdots + \alpha_v x_v$ is at least $C^* q^{-\lambda}$ in absolute value for all $x_1, \ldots, x_v \in B_1^* \cup B_2^*$ not all coming from $B_1^*$ and not all coming from $B_2^*$.
\begin{mypro}\label{simul-mainprop}
There exists $C^*$ depending only on $\alpha$ and on the local field $K$ with the following property. For any ball $B$ with $q^{-\lambda} := \text{radius}(B)$, there exist $B_1^*, B_2^* \subset B$ such that $B_1^*, B_2^*$ are balls of radius $C^*q^{-\lambda}$ and such that $\alpha_1 x_1 + \cdots + \alpha_v x_v$ has absolute value at least $C^* q^{-\lambda}$ for $x_1, \cdots, x_v \in B_1^* \cup B_2^*$ provided that not all of the $x_1, \cdots, x_v$ are in $B_1^*$, and not all of the $x_1, \ldots, x_v$ are in $B_2^*$.
\end{mypro}
We are now ready to prove Theorem \ref{simul-thm}.
\begin{proof}[Proof of Theorem \ref{simul-thm}]
This construction closely follows the one appearing in \cite{FP16}. 

In order to construct the set $E$, we begin by selecting a ball $E_0$ of radius $q^{-\lambda_0}$ where $\lambda_0$ is chosen sufficiently large so that $C q^{-\lambda_0} v < (C^*)^3$, where $C^*$ is the constant from Proposition \ref{simul-mainprop}.  

We will construct the set $E$ by applying Proposition \ref{simul-mainprop} inductively. $E_j$ will always be a union of $2^j$ balls of radius $(C^*)^j q^{- \lambda_0}$ with the property that $|\alpha \cdot (x_1, \ldots, x_v)| \geq (C^*)^j q^{- \lambda_0}$ for $x_1, \ldots, x_v \in E_j$ unless $x_1, \ldots, x_v$ all belong to the same ball in $E_j$. To construct $E_{j+1}$, we apply Proposition \ref{simul-mainprop} to each of the balls that constitute $E_j$. The result will be a collection $E_{j+1}$ of $2^j$ balls of radius $(C^*)^{j+1}q^{- \lambda_0}$ with the property that $|\alpha \cdot (x_1, \ldots, x_v)| \geq (C^*)^{j+1}q^{- \lambda_0}$ unless $x_1, \ldots, x_v$ belong to the same ball in $E_{j+1}$. Let $E = \bigcap_{j=0}^{\infty} E_j$. It follows from a routine computation that the Cantor set $E$ has positive Hausdorff dimension.

It remains to be seen that $f(x_1, \ldots, x_v)$ is nonzero for $x_1, \ldots, x_v \in E$ not identical. To see this, notice that there exists a minimal ball $B$ such that $x_1, \ldots, x_v \in B$. Say that this ball $B$ is a basic ball of $E_{j-1}$. Let $B_1$ and $B_2$ be the children of this ball in the construction. Then we have $|\alpha \cdot (x_1, \ldots, x_v)| \geq (C^*)^j q^{- \lambda_0}$.  Since $x_1, \ldots, x_v \in B_1 \cup B_2 \subset B$, we have that 
\[(|x_2 - x_1|^2 + \cdots + |x_v - x_1|^2) \leq v (C^*)^{2j - 2} q^{- 2 \lambda_0}.\]
 Therefore we have that 
\[|f(x_1, \ldots, x_v)| \geq (C^*)^j q^{- \lambda_0} > 0\]
 as long as 
\[C v (C^*)^{2j - 2} q^{-2 \lambda_0} < (C^*)^j q^{\lambda_0}\] 
by the ultrametric inequality. Noting that $C q^{- \lambda_0} v < (C^*)^3$, we get $Cv (C^*)^{2j - 2} q^{- 2 \lambda_0} \leq (C^*)^{2j +1} q^{- \lambda_0} \leq C^*{j} q^{- \lambda_0}$ as desired.
\end{proof}
\bibliographystyle{plainurl}
\bibliography{Local_Configurations}
\end{document}